\documentstyle[epsfig]{article}
\newtheorem{theorem}{Theorem}[section]
\newtheorem{lemma}{Lemma}[section]

\newtheorem{remark}{Remark}[section]
\newtheorem{remarks}{Remark}[section]
\newtheorem{definition}{Definition}[section]
\newcommand{\be}{\begin{equation}}
\newcommand{\ee}{\end{equation}}

\newcommand{\teta}{\theta}

\newcommand{\om}{\omega}

\newcommand{\ov}{\overline}

\newcommand{\wtilde}{\widetilde}

\renewcommand{\d }{\delta }

\newcommand{\g }{\gamma}

\newcommand{\vphi}{\varphi }

\begin{document}

\title{{\bf Optimal stability and instability results
for a class of nearly integrable Hamiltonian systems}}

%{Massimiliano Berti S.I.S.S.A., Via Beirut 2-4,
%34014, Trieste, Italy ({\tt berti@sissa.it}) } and Philippe Bolle
%D\'epartement de math\'ematiques, Universit\'e
%d'Avignon, 33, rue Louis Pasteur, 84000 Avignon, France
%({\tt philippe.bolle@univ-avignon.fr})}}
\author{Massimiliano Berti, Luca Biasco and Philippe Bolle}
\date{}
\maketitle

{\bf Abstract:}
We consider a  nearly integrable, non-isochronous,  
a-priori unstable Hamiltonian system
with a (trigonometric polynomial)  
$O(\mu)$-perturbation which does not preserve the unperturbed tori. 
We prove the existence of Arnold diffusion with diffusion time 
$ T_d = O((1/ \mu) \log (1/ \mu ))$ by a variational method  
which does not require the existence
of ``transition chains of tori'' provided by KAM theory.
We also prove that our estimate of the diffusion time $T_d $  
is optimal as a consequence of a general stability result
proved via classical perturbation theory.
\footnote{Supported by M.U.R.S.T. Variational Methods and Nonlinear
Differential Equations.}
\\[2mm]
Keywords: Arnold diffusion, variational methods,
shadowing theorem, perturbation theory, nonlinear functional analysis
\\[1mm]
AMS subject classification: 37J40, 37J45.
\\[2mm]
\indent
{\bf Riassunto:} {\it Risultati ottimali di stabilit\`a 
e di instabilit\`a per una classe di sistemi Hamiltoniani quasi-integrabili}.
Consideriamo sistemi Hamiltoniani quasi-integrabili, non-isocroni,
a-priori instabili soggetti ad una
pertubazione $O(\mu)$ (un polinomio trigonometrico) che non 
preserva i tori imperturbati. Mediante un metodo variazionale
dimostriamo l'esistenza di orbite di diffusione con tempo di diffusione 
$ T_d = O((1/ \mu) \log (1/ \mu ))$. Il nostro approccio
non richiede l'esistenza di ``catene di tori KAM di transizione''.
Proviamo inoltre l'ottimalit\`a
della nostra stima sul tempo di diffusione, in conseguenza di un
risultato generale di stabilit\`a dimostrato mediante la teoria
classica delle perturbazioni.  

\section{Introduction}

We outline in this Note some recent results on Arnold's diffusion 
obtained in \cite{BBB} where we refer for complete
proofs. We consider nearly integrable
{\it non-isochronous} Hamiltonian systems described by
\be\label{eq:Hamg}
{\cal H}_\mu =  \frac{I^2}{2} + \frac{p^2}{2} + (\cos q - 1) +
\mu f(I, p, \vphi, q , t),
\ee
where $( \vphi, q , t) \in {\bf T}^d \times {\bf T}^1 \times {\bf T}^1 $
are the angle variables, $ ( I, p ) \in {\bf R}^d \times {\bf R}^1$
are the action variables and $ \mu \geq 0 $ is a small real  
parameter. 
The Hamiltonian system associated with ${\cal H}_{\mu}$ writes
$$
\dot{\vphi}= I + \mu \partial_I f,  \quad
\dot{I}= -\mu \partial_{\vphi} f, \quad
\dot{q}= p + \mu \partial_p f,  \quad  
\dot{p} = \sin q - \mu \partial_q f.  \eqno({\cal S}_\mu)
$$
The perturbation $f$ is assumed, as in \cite{CG}, 
to be a trigonometric polynomial
of order $ N $ in $ \vphi $ and $ t $, namely
\be\label{trigpert}
f(I, p, \vphi, q, t) = \sum_{ |(n,l)| \leq N} f_{n,l} (I, p, q) 
exp^{ {\rm i}(n \vphi  + l t)}.
\ee
$({\cal S}_{\mu})$ describes a system of $ d $ ``rotators'' 
weakly coupled with a pendulum through a small periodically time dependent
perturbation term.
The unperturbed Hamiltonian system $ ({ \cal S }_0) $ is completely integrable and,
in particular, the energy $ I_i^2/ 2 $ of each rotator is a
constant of the  motion.
The problem of {\it Arnold diffusion} in this context is whether,
for $ \mu \neq 0 $, there exist motions whose net effect is to transfer
$O(1)$-energy among the rotators. A natural complementary question regards 
the time of stability (or instability) 
for the perturbed system: what is the minimal time to 
produce an $O(1)$-exchange of energy, if any takes place, 
among the rotators?
\\[1mm]
\indent
The mechanism proposed in \cite{Arn} to prove the existence 
of Arnold diffusion and thereafter become classical, is the following one.
The unperturbed Hamiltonian system $ ({\cal S}_0) $ admits
a continuous family of $d$-dimensional 
partially hyperbolic invariant tori
$ {\cal T}_\om = \{ \vphi  \in  {\bf T}^d,
I = \om, \ q = p = 0 \}$
possessing stable and unstable manifolds
$ W^s_0 ( {\cal T}_\om) = W^u_0 ( {\cal T}_\om ) =$
$ \{ \vphi \in  {\bf T}^d,
I = \om, \ (p^2/2) + (\cos q - 1) = 0 \}$.
The method used in \cite{Arn} to produce unstable orbits 
relies on the construction, for $\mu \neq 0, $ of  {\it ``transition chains''} 
of perturbed partially hyperbolic  
tori $ {\cal T}_\om^\mu $ close to $ {\cal T}_\om $
connected one to another by heteroclinic orbits. Therefore in general 
the first step 
 is  to prove the persistence of such hyperbolic tori 
$ {\cal T}_\om^\mu $ for $ \mu \neq 0 $ small enough, and to show that 
its perturbed stable and unstable manifolds $ W^s_\mu ({\cal T}_\om^\mu )$
and $ W^u_\mu ({\cal T}_\om^\mu) $ split and intersect transversally
(``splitting problem''). 
The second step is to find a transition chain of perturbed tori:
this is a difficult task since, for general non-isochronous systems,  
the surviving perturbed tori $ {\cal T}^\mu_\om $ are separated 
by the gaps appearing in KAM constructions.
Two perturbed invariant tori $ {\cal T}_\om^\mu $ and $ {\cal T}_{\om'}^\mu $
could be too distant one from
the other, forbidding the existence of an heteroclinic intersection
between $ W^u_\mu ({\cal T}_\om^\mu )$ and $ W^s_\mu ({\cal T}_{\om'}^\mu )$:
this is the famous {\it ``gap problem''}. 
In \cite{Arn} this difficulty is bypassed by the peculiar 
choice of the perturbation $ f ( I, \vphi, p, q,t ) = 
(\cos q - 1 ) f ( \vphi, t ) $, 
whose gradient vanishes on the unperturbed tori $ {\cal T}_\om $, leaving 
them {\it all} invariant also for $ \mu \neq 0 $.  
The final step is to prove, by a ``shadowing argument'', the existence of 
a true diffusion orbit, close to a given transition chain 
of tori, for which the action variables $I$ 
undergo a drift of $O(1)$
in a certain time $ T_d $ called the {\it diffusion time}.

The first paper proving Arnold diffusion in presence of
perturbations not preserving the unperturbed tori was \cite{CG}. 
Extending Arnold's analysis, it is proved in \cite{CG} 
that, if the perturbation is 
a trigonometric polynomial in the angles $ \vphi $, % and the time $t$, 
then, in some regions of phase space, the ``density'' of 
perturbed invariant tori is  high enough for the 
construction of a transition chain. 

Regarding the shadowing problem, geometrical method, see e.g.
\cite{CG}, \cite{M}, \cite{C1}, 
\cite{CrG},  and variational ones, see e.g. \cite{BCV}, 
have been applied, in the last years, 
in order to prove the existence of diffusion orbits shadowing 
a given transition chain of tori and to estimate the diffusion time.
We also quote the important paper \cite{Bs}  which, even if 
dealing only with the Arnold's model perturbation, 
has obtained, using variational methods, very good 
time diffusion estimates 
and has introduced new ideas for studying the shadowing problem.
For isochronous systems new variational results
concerning the shadowing and the splitting problem have been obtained in
\cite{BBN}, \cite{BB1} and \cite{BB3}. 
\\[2mm]
\indent
In this Note we describe an {\it alternative mechanism},
proposed in \cite{BBB}, to produce diffusion orbits. This 
method is not based 
on the existence of transition chains of tori, namely 
it avoids the KAM construction of the perturbed hyperbolic tori, 
proving directly the existence of 
a drifting orbit as a local minimum of an action functional,
see Theorem \ref{thm:main}.
At the same time this variational approach 
achieves the optimal diffusion time $ T_d = O(( 1 / \mu) \ln ( 1/ \mu ))$,
see (\ref{totaltime}). 
We also prove that our time diffusion estimate 
is the optimal one as a consequence of a general stability result,
Theorem \ref{thm:Nec}, proved via classical perturbation theory. 
As in \cite{CG} our diffusion orbit will not connect any two arbitrary
frequencies of the action space, even if we manage to connect
more frequencies than in \cite{CG}, proving the drift also
in some regions of phase space  
where transition chains might not exist.
Clearly if the perturbation is chosen as in  Arnold's example 
we can drift in all phase the  space with no restriction, see Theorem 
(\ref{thm:Arn}).
\\[1mm]
\indent
Actually our variational shadowing technique is not 
restricted to the a-priori unstable case, but would
allow, in the same spirit of \cite{BBN}, \cite{BB1} and \cite{BB3}, 
once a ``splitting property'' is somehow proved,
to get diffusion orbits with the best diffusion time
(in terms of some measure of the splitting).
\\[1mm]
\indent
In conclusion the results and the method described in this Note
constitute a further step in a research line,
started in \cite{BBN}-\cite{BB1} and \cite{BB3}, 
whose aim is to find new mechanisms for proving Arnold diffusion.
We expect that these variational methods %developed in \cite{BBB} 
could be suitably refined in order to prove the existence of 
drifting orbits in the whole phase space, and also 
for generic analytic perturbations.
Another possible application of these methods could regard 
infinite dimensional Hamiltonian systems where the existence of 
``transition chains of infinite dimensional hyperbolic tori'' 
is far for being proved.
\\[2mm]
\indent
{\bf Acknowledgments:} 
Part of this paper was written when the third
author was visiting SISSA in Trieste. 
He thanks Prof. A. Ambrosetti and  the staff of SISSA 
for their kind hospitality.

\section{Main results}

For simplicity, even if not really necessary, when proving 
the existence of diffusion orbits, we assume $f$ to be 
a purely spatial perturbation, namely 
$f(\vphi, q, t) = \sum_{ |(n,l)| \leq N} f_{n,l} (q) 
exp ({\rm i}(n \vphi  + l t)) $. 
The functions $ f_{n,l} $ are assumed to be  smooth.

Let us  define the ``resonant web'' $ {\cal D}_N $, formed by the frequencies
$ \om $ ``resonant with the perturbation'',  
\be\label{eq:resonantweb}
{\cal D}_N := \Big\{ \om \in {\bf R}^d \ \Big| \ \exists (n, l)
\in {\bf Z}^{d+1} \ {\rm s.t.} \
0 < | ( n, l ) | \leq N \ {\rm and } \ \om n + l = 0 \Big\} = 
\cup_{0< |(n,l)| \leq N} E_{n,l}
\ee
where $ E_{n,l} := \{ \om \in {\bf R}^d \ | \ \om n + l = 0 \} $. 
Let us also consider the Poincar\'e-Melnikov primitive 
$$
\Gamma (\om, \vphi_0 , \teta_0 ) := - \int_{\bf R}
\Big[ f ( \omega t + \vphi_0, q_0 (t), t + \teta_0) -  
f ( \omega t + \vphi_0, 0 , t + \teta_0 ) \Big] \ dt, 
$$ 
where $ q_0 (t) = 4 \ {\rm arctan} (\exp t ) $ is
the unperturbed separatrix of the pendulum satisfying $ q_0 (0) = \pi $.

The next theorem states that, for any 
connected component $ {\cal C} \subset {\cal D}_N^c $,
$ \om_I, \om_F \in {\cal C} $, 
there exists a solution of $ ({\cal S}_\mu) $
connecting a $O(\mu)$-neighborhood of $ {\cal T}_{\om_I} $ to a 
$O(\mu)$-neighborhood of ${\cal T}_{\om_F}$, 
in a time-interval of length 
$ T_d = O((1 / \mu) | \log \mu |) $.

\begin{theorem}\label{thm:main}
Let ${\cal C}$ be a connected component of $ {\cal D}_N^c $, 
$\om_I, \om_F \in {\cal C} $ and let 
$ \gamma : [0,L] \to {\cal C} $ be a smooth embedding such that
$ \g (0) = \om_I $ and $\gamma (L) = \om_F $.
Assume that $\forall \om := \gamma (s) $ $(s \in [0,L])$,  
$ \Gamma(\om, \cdot, \cdot) $ possesses a 
non-degenerate minimum $ ( \vphi_0^\om, \teta_0^\om ) $. 
Then $\forall \eta > 0 $ there exists $ \mu_0 := \mu_0 (\g, \eta) > 0 $
and $C:= C(\gamma) >0 $ such that $ \forall 0 < \mu \leq \mu_0 $ 
there exists a solution 
$ (\vphi_\mu (t), q_\mu (t),$ $ I_\mu (t) , p_\mu (t))$ of  
$ ({\cal S}_\mu) $ and two instants 
$ {\tau}_1 < {\tau}_2 $ such that $ I_\mu ( \tau_1 ) = \om_I + O(\mu) $,  
$ I_\mu ( \tau_2 ) = \om_F + O(\mu) $, 
\be\label{totaltime}
|{\tau}_2 - {\tau}_1| \leq  
\frac{C}{ \mu } | \log \mu |
\ee
and ${\rm dist} (I_\mu (t), \gamma ) < \eta $ 
for all $\tau_1 \leq t \leq \tau_2 $. 
%Finally, the previous result holds 
%for any smooth perturbation $ f + {\wtilde f}$ with 
%${\wtilde f} = O( \mu^a ) $, $ a > 0$.
\end{theorem} 

We can also build diffusion orbits approaching
the boundaries of $ {\cal D}_N $ at distances as 
small as a certain power of $ \mu $: see for a precise statement 
Theorem 6.1 of \cite{BBB}.
Theorem \ref{thm:main} improves the corresponding result in \cite{CG}
which enables to connect any two frequencies $ \om_I $ and $ \om_F $ 
belonging to the same connected component
${\cal C} \subset {\cal D}_{N_1}^c $ for $ N_1 = 14 d N $ and 
with dist$( \{ \om_I, \om_F \}, {\cal D}_{N_1}) = O( 1 )$. 
%Unlikely we are able to connect any two frequencies $ \om_I $ and $ \om_F $
%inside the same connected component 
%${\cal C} $ of $ {\cal D}_N^c $ which is strictly larger than 
%any connected components of $ {\cal D}_{N_1}^c $.
As already said such restriction of \cite{CG}  
arises because transition chains might not exist in the whole 
$ {\cal C} \subset {\cal D}_N^c $. 
%see also remark (\ref{remtranl})-($iv$).  

 Theorem \ref{thm:main} 
improves also the known time-estimates on the diffusion time.
The first estimate on the diffusion time obtained
by geometrical method in \cite{CG} is $ T_d = O( \exp{(1/\mu^2)} )$.
%and it has been improved in \cite{G1} to be $ T_d = O( \exp{(1/\mu)} )$.
In \cite{M}-\cite{C1}-\cite{CrG}, 
still by geometrical methods, and in \cite{BCV}, 
by means of Mather's theory, the diffusion time 
has been proved to be just polynomially long in the splitting $ \mu $ 
(the splitting angles between the perturbed stable and unstable manifolds
$ W^{s, u}_\mu ({\cal T}_\om^\mu ) $
at a homoclinic point are, by classical Poincar\'e-Melnikov theory, 
$ O(\mu )$).
We note that the variational method proposed by Bessi in \cite{Bs} had 
already given, even if in the case of perturbations preserving 
all the  unperturbed tori, the time diffusion estimate $ T_d = O(1/ \mu^2) $. 
For isochronous systems the estimate on the 
diffusion time  $ T_d = O((1/ \mu) |\ln \mu |) $ has already been obtained in 
\cite{BBN}-\cite{BB1}. Very recently, in \cite{CrG}, the diffusion time 
has been estimated as $ T_d = O((1/\mu) |\log \mu |)$ by a method 
which uses ``hyperbolic periodic orbits''; 
however the result of \cite{CrG} is of local nature: 
the previous estimate holds only for 
diffusion orbits shadowing a transition chain close to some torus run with
diophantine flow. 

Our next statement (a stability result) concludes this  quest of the minimal 
diffusion time $ T_d $ : it proves the 
%Clearly Theorems \ref{thm:Nec} and \ref{thm:main} 
optimality of our estimate $ T_d = O( (1 / \mu) |\log \mu |)$.
%(\ref{eq:stab}) and 
%(\ref{totaltime}) for the diffusion time.
%is exactly the optimal's one. 
%Theorem \ref{thm:Nec},
%proved in section \ref{sec:opt},
%gives a lower bound on the minimal time for 
%the rotators to exchange, {\it if ever possible}, $O(1)$ energy. 
%Indeed, our next instability result 
%is the exact complement of Theorem \ref{thm:Nec}:
%we prove that Arnold diffusion
%actually takes place, in some regions of phase space, 
%in  the minimal diffusion time 
%$T_d = O((1/ \mu) |\ln \mu |$). 

\begin{theorem}\label{thm:Nec}
Let $ f ( I, \vphi, p, q, t) $ be as in (\ref{trigpert}), 
where the $f_{n,l}$ ($|(n,l)|\leq N$) are analytic functions. Then $ \forall \kappa, 
\ov{r}, \wtilde{r} > 0 $ there exist 
$ \kappa_0, \mu_1 > 0 $ such that $ \forall \, 0 < \mu \leq \mu_1 $,
any solution $ ( I(t), \vphi(t), q(t), p(t))$ of $({\cal S}_\mu)$
with $| I(0) | \leq \ov{r} $ and $|p(0)| \leq \wtilde{r} $
satisfies
\be\label{eq:stab}
|I(t)-I(0)| \leq  \kappa, \ \qquad \ \forall\ 
|t| \leq \frac{\kappa_0}{\mu} \ln \frac{1}{\mu}.
\ee
\end{theorem}

Actually the proof of Theorem \ref{thm:Nec} contains much more
information: in particular the stability time (\ref{eq:stab}) 
is sharp only for orbits lying close to the separatrices. On the other hand
the orbits lying far away from the separatrices are much more
stable, namely exponentially stable in time according to Nekhoroshev 
type time estimates, see (\ref{stabexpi}).
Indeed the diffusion orbit of Theorem \ref{thm:main} is found 
(as a local minimum of an action functional) 
close to some pseudo-diffusion orbit whose $(q,p)$ variables 
stay close to the separatrices of the pendulum
(turning $O(1/ \mu )$ times around them).
\\[1mm]
\indent
As a byproduct of the techniques developed in this paper 
we have the following result concerning  ``Arnold's example'' 
\cite{Arn}

\begin{theorem}\label{thm:Arn}
Let $ f(\vphi, q, t) := (1 - \cos q) {\wtilde f} (\vphi, t) $. Assume that 
%$ \forall \om_I, \om_F \in {\bf R}^d $ we consider 
for some smooth embedding $ \gamma : [0,L] \to {\bf R}^d $, 
with $\g (0) = \om_I $ and $\gamma (L) = \om_F $, 
$\forall \om := \gamma (s) $ $(s \in [0,L])$,  
$ \Gamma(\om, \cdot, \cdot) $ possesses a 
non-degenerate minimum $ ( \vphi_0^\om, \teta_0^\om ) $. 
Then $ \forall \eta > 0 $ there exists $ \mu_0 := \mu_0 ( \g , \eta ) > 0 $,
and $ C:= C(\gamma) >0 $ such that $\forall 0 < \mu \leq \mu_0 $ 
there exists a heteroclinic orbit ($\eta$-close to $ \gamma $) 
connecting the invariant tori $ {\cal T}_{\om_I}$ and $ {\cal T}_{\om_F} $.
Moreover the diffusion time $ T_d $ needed
to go from a $\mu$-neighbourhood of
${\cal T}_{\om_I}$ to a $\mu$-neighbourhood of
${\cal T}_{\om_F}$ is bounded by
$T_d \leq  (C / \mu)  | \log \mu |$.
\end{theorem} 

The method of proof of Theorem \ref{thm:main} (and Theorem \ref{thm:Arn}) 
relies on a finite dimensional reduction 
of Lyapunov-Schmidt type, variational in nature, 
introduced in \cite{AB} and later extended in \cite{BBN},\cite{BB1} 
and \cite{BB3} to the problem of Arnold diffusion.
The diffusion orbit of Theorem \ref{thm:main} 
(and Theorem \ref{thm:Arn}) is found 
as a local minimum of the action functional 
close to some pseudo-diffusion orbit whose $(q,p)$ variables move 
along the separatrices of the pendulum. The pseudo-diffusion orbits, 
constructed by the Implicit Function Theorem, are
true solutions of $ ({\cal S}_\mu) $ 
except possibly at some instants $ \teta_i $, for $ i = 1, \ldots, k $, 
when they are glued 
continuously at the section $\{ q = \pi , \ {\rm mod} \ 2 \pi {\bf Z} \}$ 
but the speeds $({\dot q}_\mu (\teta_i), 
{\dot \vphi}_\mu (\teta_i) ) = (p_\mu (\teta_i), 
I_\mu (\teta_i) ) $ may possibly have a jump, see lemma 3.1 of \cite{BBB}.
The time interval $ T_s = \teta_{i+1} - \teta_i $ is heuristically the
time required to perform a single transition during which the rotators 
can exchange $ O(\mu)$-energy, i.e. the action variables vary of $O(\mu)$.
During each transition we can exchange only $O(\mu )$-energy 
since the splitting is $O(\mu )$-large. Hence 
%Since during each transition we can exchange only $ O ( \mu ) $ energy,
in order to exchange $O(1)$ energy  the number of transitions required 
will be $k = O(1/ {\rm splitting}) = O(1/ \mu)$.
\\[1mm]
\indent
We underline that the question of finding the optimal time and 
the mechanism for which we can avoid the construction of true transition chains
of tori are deeply connected. Indeed the main 
reason for which our drifting technique  avoids the construction 
of KAM tori is the following one: if
the time to perform a simple transition $ T_s $
is, say, just $ T_s = O( |\ln \mu | ) $  then, on such short time intervals, 
it is easy to approximate the action functional with the unperturbed solutions 
living on the stable and unstable manifolds of the unperturbed tori 
$ W^s ({\cal T}_\om ) = W^u({\cal T}_\om ) = 
\{ (\vphi , \om , q,p) \ | \ 
  p^2 / 2 + ( \cos q - 1) = 0 \}$, see lemma 5.6 of \cite{BBB}. 
In this way we do not need 
to construct the true hyperbolic tori $ {\cal T}_\om^\mu $
close to $ {\cal T}_\om $
(actually for our approximation we only need the time for a 
single transition to be $ T_s << 1/ \mu $).
\\[1mm]
\indent
The fact that it is possible to perform a single transition in a 
very short time interval like $ |\ln \mu |$ is not obvious at all. 
In \cite{Bs} the time to perform a single transition,
in the example of Arnold, is $O(1 / \mu )$. 
This time separation arises in order to ensure that the 
variations of the action functional of the rotators are small compared with 
the (positive definite) second derivative  
of the Poincar\'e-Melnikov primitive at its minimum point.
Unfortunately this time is too long  to make directly our 
approximations of the action functional. 
The key observation that enables us to perform 
a single transition in a very short time interval  
concerns the behaviour the ``gradient flow'' 
of the unperturbed action functional of the rotators. 
In section 6 of \cite{BBB} it is shown that 
the variations of the action of the rotators are small, 
even on time intervals $T_s << 1/ \mu $, and do not ``destroy'' the minimum
of the Poincar\'e-Melnikov primitive. 
\\[1mm]
\indent
When trying to build a pseudo-diffusion orbit which 
performs single transitions in very short time intervals 
we encounter another difficulty linked with the ergodization time. 
The time to perform a single transition $T_s$ must be sufficiently long 
 to settle, at each instant $ \teta_i $, the projection of the 
pseudo-orbit on the torus sufficiently close
to the minimum of the Poincar\'e-Melnikov function, i.e. the homoclinic point
(in our method 
it is sufficient to arrive just $O(1)$-close, independently of $ \mu$,
to the homoclinic point). This necessary request creates some
difficulty since our pseudo-diffusion orbit 
may arrive $O(\mu)$-close in the action space 
to resonant hyperplanes of frequencies
whose linear flow does not provide a dense enough net of the torus. 
The way in which this problem is overcome is discussed in 
\cite{BBB}: 
we observe a phenomenon of ``stabilization close to resonances''
which forces the time $ T_s $  for some single transitions  
to increase . 
Anyway the total time required to cross these
(finite number of) resonances is still 
$ T_d = O( (1/ \mu) \log (1 / \mu)) $.
%see the proof of Theorem \ref{thm:main}. 
This  discussion  enables us to prove optimal 
fast-Arnold diffusion in large regions of phase space and allows to improve 
the local diffusion results of \cite{CrG}.
\\[1mm]
\indent
As explained before we need, 
in order to prove Theorems \ref{thm:main} and \ref{thm:Arn},
some estimates on the ``ergodization time 
of the torus'' for linear flows possibly resonant but 
only at a ``sufficiently high order''. 
The following result (Lemma \ref{lem:ergo2}) gives an answer which may be
 of independent interest.
Let $\Gamma$ be a lattice of ${\bf R}^l$, {\it i.e.} a discrete
subgroup of ${\bf R}^l$ such that $ {\bf R}^l / \Gamma$ has finite volume.
$\forall \Omega \in {\bf R}^l $ the ``ergodization time 
required by the flow $ \{ \Omega t \} $ 
to fill the torus $ {\bf R}^l / \Gamma $ within $ \d $'' is defined as
$ T ( \Gamma, \Omega, \delta ) := \inf \{ t\in {\bf R}_+ \ | 
\ \forall x \in {\bf R}^l \  d(x, [0,t] \Omega + \Gamma) \leq \delta 
\}$ ( with inf$E = +\infty$ if $E=\emptyset$). 
For $ R > 0 $ define also $
\Gamma^*_R = \{ p \in {\bf R}^l \ | \ 0 < |p| \leq R \ , \ 
\forall \gamma \in \Gamma \ 
p \cdot \gamma \in {\bf Z} \}$ and 
$ \alpha(\Gamma,\Omega, R)= 
\inf \{ |p \cdot \Omega | \ | \  p \in \Gamma^*_R \}.$
The following result holds:
\begin{lemma}\label{lem:ergo2}
$\forall l \in {\bf N}$, $ \exists a_l > 0 $ 
such that, for all lattice $ \Gamma \subset {\bf R}^l $, 
$ \forall \Omega \in {\bf R}^l$,  
$\forall \delta >0$,  
$ T ( \Gamma, \Omega, \delta) \leq (\alpha (\Gamma,\Omega, a_l/\delta))^{-1}$.
Moreover $ T (\Gamma,  \Omega, \delta ) \geq ( 1 / 4 ) 
\alpha( \Gamma, \Omega, 1/4\delta)^{-1} $.
\end{lemma}

A straighforward consequence of Lemma \ref{lem:ergo2} is, for example,
Theorem D of \cite{BGW}:
if $ \Omega $ is a $C-\tau$ diophantine vector, i.e. 
there exist $ C > 0 $ and $ \tau \geq l-1 $ such that 
$| k \cdot \Omega| \geq C / |k|^{\tau}$, 
$ \forall k \in {\bf Z}^l \setminus \{ 0 \}$, then 
$ T(\Omega, \delta) \leq a_l^{\tau}/(C\delta^{\tau})$
(where $T(\Omega, \delta) := T({\bf Z}^l , \Omega, \delta)$). 
Also Theorem B of \cite{BGW} is an easy consequence of Lemma \ref{lem:ergo2}.
\\[2mm]
\indent
We conclude this Note discussing briefly the 
proof of  Theorem \ref{thm:Nec}, see section 7 of \cite{BBB}.
First we prove stability in the region ``far from the 
separatrices of the pendulum'' 
$ {\cal E}_1 := \{ (I, \vphi, q,p) \ | \ | E(q,p) | \geq \mu^c \} $, 
where $ E(q,p) := p^2/2 + (\cos q - 1) $  and $ c $ is a suitable
positive constant.
In $ { \cal E }_1 $ we can write
Hamiltonian $ {\cal H}_\mu $ in action-angle variables 
$ (I,\vphi, Q, P, t)$ where $Q:= Q(q,p)$ and  $P:= P(q,p)$ 
are the action-angle variables of the standard pendulum 
$ E(q,p) $, i.e $ E(q(Q,P),p(Q,P)) := K(P) $. 
In these variables the new Hamiltonian writes 
$ {\cal H}_1 :=  I^2 / 2 + K(P) + \mu 
f_1( \vphi, Q , t, P, I) $, and, by a result of \cite{BC}, it is 
analytic with an analyticity radius $ r_1 \approx \mu^c $
(when $\mu $ goes to zero, the region  ${\cal E}_1 $ 
approximate closer and closer to the separatrices and the 
analyticity estimate deteriorates).
It turns out that $ {\cal H}_1 $ is {\sl steep} (actually
for $E $ positive it is even  {\sl quasi-convex}, see \cite{Po}) 
and then, for $c>0$ small enough, we can 
apply the Nekhoroshev Theorem as proved in 
\cite{N}. In this way we obtain  
exponential stability in the whole region ${ \cal E }_1 $, i.e.
\be\label{stabexpi}
| I(t)- I (0) | \leq const. \mu^b   \ \qquad \ \forall \ 
| t | \leq T := const \frac{1}{\mu}\exp \Big( \frac{1}{\mu} \Big)^a 
\ee
for two constants $ a,b > 0 $. 
Finally we study the behaviour of an orbit close to the separatrices of 
the pendulum, namely in the  region ${ \cal E }_1^c $.
Roughly speaking, such an orbit will spend alternatively a time 
$T_S = O( |\ln \mu |) $ into a small O(1)-neighborhood of 
$(p,q) = (0,0)$ and a time $ T_F = O(1) $ outside. In this second case 
we directly obtain
$ \Delta I := \int_{s}^{s+T_F} {\dot I}_\mu = O(\mu T_F) = O( \mu )$. 
In \cite{BBB} it is proved that  
$ \Delta I = O ( \mu ) $ also in the first case. 
This result is obtained performing one step of classical 
perturbation theory (as in \cite{CG}) after writing the pendulum
in hyperbolic variables in a small neighborhood of the origin.

\noindent
{\it Massimiliano Berti and Luca Biasco, S.I.S.S.A., Via Beirut 2-4,
34014, Trieste, Italy, berti@sissa.it}.
\\[2mm]
{\it Philippe Bolle,
D\'epartement de math\'ematiques, Universit\'e
d'Avignon, 33, rue Louis Pasteur, 84000 Avignon, France,
philippe.bolle@univ-avignon.fr}

\end{document}